\documentclass[11pt,a4paper]{article}

\usepackage[english]{babel}
\usepackage[leqno]{amsmath}
\usepackage{amssymb,amsthm}

\def\hyp{\mathcal H}
\def\hypbar{\overline{\mathcal H}}
\def\calo{\mathcal O}
\def\calz{\mathcal Z}
\def\Pic{{\rm Pic}}
\def\QQ{\mathbb Q}

\def\CC{\mathbb C}
\def\PP{\mathbb P}

\begin{document}

\newtheorem{thm}{Theorem}
\newtheorem{prop}[thm]{Proposition}
\newtheorem{lemma}[thm]{Lemma}
\newtheorem{cor}[thm]{Corollary}

\title{The Picard group of the moduli stack\\
of stable hyperelliptic curves}
\author{Maurizio Cornalba
\footnote{partially supported by PRIN 2005 \textit{Spazi di moduli e teorie di Lie} and FAR 2004 (Pavia) \textit{Variet\`a algebriche, calcolo algebrico, grafi orientati e topologici}.}\\
{\small Dipartimento di Matematica ``F.~Casorati''}\\
{\small Universit\`a di Pavia}\\
{\small Via Ferrata 1, 27100 Pavia, Italia}
}
\date{}

\maketitle

In a recent paper \cite{arsivi}, Arsie and Vistoli have shown, as a byproduct of their study of moduli of cyclic covers of projective spaces, that the Picard group of the moduli stack $\hyp_g$ of smooth hyperelliptic curves of genus $g\ge 2$ is finite cyclic, and that its order is $8g+4$ for odd $g$, and $4g+2$ for even $g$, over any field of characteristic not dividing $2g+2$. In another recent paper \cite{gorvi}, Gorchinskiy and Viviani have given a geometric construction of generators for the Picard groups in question. On the other hand, it was shown in \cite{corhar} that, in characteristic zero, the identity
\begin{equation}
\label{oldformula}
(8g+4)\lambda=g\xi_{irr}+2\sum_{i=1}^{\lfloor\frac{g-1}{2}\rfloor}
(i+1)(g-i)\xi_i+4\sum_{j=1}^{\lfloor\frac{g}{2}\rfloor}
j(g-j)\delta_j
\end{equation}
holds in $\Pic(\hypbar_g)\otimes \QQ$, where $\hypbar_g$ is the closure of $\hyp_g$ inside the stack $\overline{\mathcal M}_g$ of genus $g$ stable curves. In the formula, $\lambda$ stands for the Hodge class, while $\xi_{irr},\xi_1,\dots,\xi_{\lfloor\frac{g-1}{2}\rfloor},\delta_1,\dots,\delta_{\lfloor\frac{g}{2}\rfloor}$, henceforth called \textit{boundary classes}, are the classes of the irreducible components of the complement of $\hyp_g$ in $\hypbar_g$ (see below for precise definitions). Moreover, it was proved in the same paper that the boundary classes are independent in $\Pic(\hypbar_g)\otimes \QQ$. We wish to show that, combining these results of \cite{corhar} with those of \cite{arsivi} and an idea of \cite{gorvi}, one gets almost immediately a complete description of $\Pic(\hypbar_g)$; in particular, one finds that (\ref{oldformula}) is valid already in $\Pic(\hypbar_g)$, and not just modulo torsion. One also gets an alternate - and to me somewhat simpler - proof of the result of Gorchinskiy and Viviani mentioned above. We shall work over $\CC$; however, let us mention that (\ref{oldformula}) has been proved in any characteristic by K.~Yamaki \cite{yamaki} (cf.~also \cite{kausz}), whose methods could probably be used to push everything through in all characteristics except those not covered by Arsie and Vistoli's results.

We begin by recalling some known facts about stable hyperelliptic curves and their moduli. First of all, $\hypbar_g$ is a smooth Deligne-Mumford stack of dimension $2g-1$. Its boundary, that is, the complement in $\hypbar_g$ of the dense open substack $\hyp_g$, is a divisor with normal crossings
$$
\Xi_{irr}+\Xi_1+\dots+\Xi_{\lfloor\frac{g-1}{2}\rfloor}+\Delta_1+\dots+\Delta_{\lfloor\frac{g}{2}\rfloor}
$$
(cf.~\cite{corhar}). Here the summands are the irreducible components of the boundary, which can be described as follows. Any stable hyperelliptic curve $C$ of genus $g$ has a semistable model $C'$ which can be represented as an admissible double covering (cf.~\cite{harmum}) of a stable $2g+2$-pointed curve of arithmetic genus zero; let $C'\to \Gamma$ be this covering. Any node $p$ of $\Gamma$ divides it in two pieces, one containing $j\ge 2$ marked points, and the other $2g+2-j\ge 2$ marked points; we may assume that $j\le g+1$. If $j=2i+1$ is odd, there is just one node of $C'$ lying above $p$, and we will say that it is a \textit{node of type} $\Delta_i$. Such a node divides $C'$ in two pieces of genera $i$ and $g-i$. We will say that the piece of genus $i$ is a \textit{tail of type} $\Delta_i$; of course, when $i=g/2$, the other piece is a tail of type $\Delta_i$ as well. If instead $j=2i+2$ is even, there are two nodes above $p$, and we will say that they form a \textit{pair of nodes of type} $\Xi_i$. Such a pair divides $C'$ in two pieces of genera $i$ and $g-i-1$. We will say that the piece of genus $i$ is a \textit{tail of type} $\Xi_i$; the other piece is also a tail of type $\Xi_i$ when $i=(g-1)/2$.

One passes from $C'$ to the stable model $C$ by contracting all tails of type $\Xi_0$, which are smooth rational components meeting the rest of $C'$ in just two points. The resulting nodes are said to be \textit{of type} $\Xi_{irr}$. The remaining nodes of $C'$ remain unchanged in $C$, and hence can be classified as nodes of type $\Delta_i$ or pairs of nodes of type $\Xi_i$, $i\ge 1$. Clearly, one can speak of tails of type $\Delta_i$ or $\Xi_i$, $i\ge 1$, also for $C$. We now define the divisors $\Xi_{irr}$, $\Xi_i$, and $\Delta_i$, for $i\ge 1$, as the loci of stable hyperelliptic curves possessing, respectively, a node of type $\Xi_{irr}$, $\Xi_i$, or $\Delta_i$.

A stable hyperelliptic curve $C$ comes with a hyperelliptic involution, which corresponds to the sheet interchange in the covering $C'\to \Gamma$. Under the involution, the nodes belonging to a pair of type $\Xi_i$ get interchanged, while nodes of all other types stay fixed.

One defines $\xi_{irr}, \xi_i, \delta_i$ as the classes in $\Pic(\hypbar_g)$ of the line bundles $\calo(\Xi_{irr})$, $\calo(\Xi_i)$, $\calo(\Delta_i)$, respectively. The classes $\delta_i$ are the pullbacks to $\hypbar_g$ of the classes with the same names in $\Pic(\overline{\mathcal M}_g)$. Instead, if we denote by $\delta_{irr}$ the pullback to $\hypbar_g$ of the class with the same name in $\Pic(\overline{\mathcal M}_g)$ (often also called $\delta_0$), then
$$
\delta_{irr}=\xi_{irr}+2\sum_{i=1}^{\lfloor\frac{g-1}{2}\rfloor}\xi_i
$$
Let $h$ be the order of the Hodge class in $\Pic(\hyp_g)$. By the results of Arsie and Vistoli, $h$ divides $8g+4$. On the other hand, since $h\lambda$ restricts to the class of a trivial line bundle on $\hyp_g$, it must be an integral linear combination of boundary classes. Since the boundary classes are independent, this relation must be proportional to (\ref{oldformula}). Thus the integer $(8g+4)/h$ divides both $g$ and $8g+4=4(2g+1)$. If $g$ is odd, it is prime with $4(2g+1)$, so the only possibility is that $h=8g+4$. We conclude that $\Pic(\hyp_g)$ is generated by the Hodge class. Now, if $\mu$ is an element of $\Pic(\hypbar_g)$, its restriction to $\hyp_g$ must be of the form $n\lambda$ for some integer $n$. Thus $\mu-n\lambda$ is an integral linear combination of boundary classes. This means that $\lambda$ and the boundary classes generate $\Pic(\hypbar_g)$. Again by the independence of the boundary classes, any relation between them and $\lambda$ must be a multiple of (\ref{oldformula}).

When $g$ is even, $h$ divides $4g+2$, by \cite{arsivi}. Moreover, all the coefficients of (\ref{oldformula}) are even. Arguing as in the odd genus case, we conclude that
\begin{equation}
\label{halfformula}
(4g+2)\lambda=\frac{g}{2}\xi_{irr}+\sum_{i=1}^{\lfloor\frac{g-1}{2}\rfloor}
(i+1)(g-i)\xi_i+2\sum_{j=1}^{\lfloor\frac{g}{2}\rfloor}
j(g-j)\delta_j
\end{equation}
If, in addition, $g$ is not divisible by $4$, the coefficients of (\ref{halfformula}) are relatively prime. Reasoning as in the odd genus case, one concludes that $\Pic(\hyp_g)$ is generated by the Hodge class, and that $\Pic(\hypbar_g)$ is generated by $\lambda$ and by the boundary classes, subject to the single relation (\ref{halfformula}).

We may summarize what has been proved in the following statement.
\begin{prop}\label{hodgegen}
Let $g\ge 2$ be an integer. When $g$ is not divisible by $4$, $\Pic(\hypbar_g)$ is generated by $\lambda$ and by the boundary classes. The relations between these classes are generated by (\ref{oldformula}) when $g$ is odd, and by (\ref{halfformula}) when $g$ is even. Moreover, (\ref{halfformula}) is valid for any even $g$.
\end{prop}
When $g$ is divisible by $4$, things are not as straightforward, since, as observed by Gorchinskiy and Viviani, in this case $\Pic(\hyp_g)$ is not generated by the Hodge class. The argument used in the preceding cases just shows that the order of the Hodge class is either $2g+1$ (the correct answer) or $4g+2$, but does not enable one to pin it down.
In order to handle this case, following \cite{gorvi} we introduce another natural line bundle on $\hypbar_g$, for any $g\ge 2$.
For any family $\alpha:X\to S$ of stable hyperelliptic curves of genus $g$, we let $W=W_\alpha$ be the divisor swept out by the Weierstrass points in the fibers, that is, the fixed scheme of the hyperelliptic involution minus the nodes of type $\Delta_i$ in the fibers. Clearly, away from nodes of type $\Xi_{irr}$, $W$ is a Cartier divisor, etale over $S$. Actually, $W$ is Cartier everywhere. In fact, in suitable local analytic coordinates, $X$ can be described near a node of type $\Xi_{irr}$ as the locus in $\CC^2\times S$ with equation $xy=f$, where $f$ is a function on $S$, while the hyperelliptic involution corresponds to $(x,y)\mapsto (y,x)$; hence a local equation for $W$ is $x-y$. Now look at the line bundle $\omega_\alpha^{g+1}(-(g-1)W)$ on $X$, where $\omega_\alpha$ is the relative dualizing sheaf, and observe that its restriction to any smooth fiber of $\alpha$ is trivial. To see this it suffices to consider the case when $S$ is a point, in which our claim is obviously true, since $\omega_X=\pi^*\omega_{\PP^1}(W)$ and $\calo(W)=\pi^*(\calo_{\PP^1}(g+1))$, where $\pi:X\to \PP^1$ is the hyperelliptic double covering.

To define a line bundle $\calz$ on $\hypbar_g$ we need to give, for each family $\alpha:X\to S$ as above, a line bundle $\calz_{\alpha}$ on $S$, natural under morphisms of families. Actually, it suffices to do this only when $S$ is etale over $\hypbar_g$. The idea would be to take as $\calz_{\alpha}$ the direct image of $\omega_\alpha^{g+1}(-(g-1)W)$. Unfortunately, $\omega_\alpha^{g+1}(-(g-1)W)$ is not necessarily trivial on singular fibers of $\alpha$, so in general this procedure does not yield a line bundle. To cure this, we twist  $\omega_\alpha^{g+1}(-(g-1)W)$ by a suitable divisor whose support is contained in the union of singular fibers. We let $G_i$ and $E_i$ be the divisors in $X$ swept out, respectively, by tails of type $\Delta_i$ and of type $\Xi_i$. We claim that $M=\omega_\alpha^{g+1}(-(g-1)W-\sum(2g-4i)G_i-\sum(g-2i-1)E_i)$ is trivial on every fiber of $\alpha$, smooth or singular. In fact, let $C$ be a fiber, and let $C'\to \Gamma$ be the corresponding admissible double covering. To show that the restriction of $M$ to $C$ is trivial it suffices to show that this is true for its pullback to $C'$. On the other hand, it is clear that this pullback comes from a line bundle on $\Gamma$; thus it suffices to show that $M$ has degree zero on every component of $C$ or, equivalently, on every tail of type $\Xi_i$ or $\Delta_i$ of $C$. This is immediate. For instance, let $T$ be a tail of type $\Xi_i$. Notice that the restriction of $\calo(E_i)$ to $T$ has degree $-2$, while the restriction of $\omega_C$ has degree $2i$. Therefore the restriction of $M$ to $T$ has degree equal to $2i(g+1)-(g-1)(2i+2)+2(g-2i-1)=0$. Similar considerations apply to the tails of type $\Delta_i$.
In conclusion,
\begin{align}
\label{defzeta}
\calz_{\alpha}&=\alpha_*(M)\\
&=\alpha_*\left(\omega_\alpha^{g+1}\left(
-(g-1)W-\sum(2g-4i)G_i-\sum(g-2i-1)E_i\right)\right)\notag
\end{align}
is a line bundle; as it behaves nicely under morphisms of families, this defines a line bundle $\calz$ on $\Pic(\hypbar_g)$. We define $\zeta\in \Pic(\hypbar_g)$ to be the class of $\calz$. When $g$ is odd, $g+1$, $g-1$, $2g-4i$, and $g-2i-1$ are all even, and replacing them with their halves in (\ref{defzeta}) produces another line bundle $\calz'$ on $\hypbar_g$, whose class we denote by $\zeta'$. Clearly, $\zeta=2\zeta'$.

We now have at our disposal all the necessary ingredients to state our main result, which is the following.

\begin{thm}\label{mainthm}
Let $g\ge 2$ be an integer. When $g$ is even, $\Pic(\hypbar_g)$ is generated by $\zeta$ and by the boundary classes $\xi_{irr},\xi_1,\dots,\xi_{\lfloor\frac{g-1}{2}\rfloor},\delta_1,\dots,\delta_{\lfloor\frac{g}{2}\rfloor}$, subject to the single relation
\begin{equation}\label{zetaeven}
(4g+2)\zeta=\xi_{irr}+2\sum_{i=1}^{\lfloor\frac{g-1}{2}\rfloor}
(i+1)(2i+1)\xi_i+4\sum_{j=1}^{\lfloor\frac{g}{2}\rfloor}
j(2j+1)\delta_j
\end{equation}
When $g$ is odd, $\Pic(\hypbar_g)$ is generated by $\zeta'$ and by the boundary classes, subject to the single relation
\begin{equation}\label{zetaodd}
(8g+4)\zeta'=\xi_{irr}+2\sum_{i=1}^{\lfloor\frac{g-1}{2}\rfloor}
(i+1)(2i+1)\xi_i+4\sum_{j=1}^{\lfloor\frac{g}{2}\rfloor}
j(2j+1)\delta_j
\end{equation}
In particular, $\Pic(\hypbar_g)$ is free abelian of rank $g$ for any $g\ge 2$.
\end{thm}
The essential step in establishing the theorem is the following result, whose proof will be given later.

\begin{lemma}\label{mainlemma}
Relation (\ref{zetaeven}) holds in $\Pic(\hypbar_g)\otimes\QQ$.
\end{lemma}
Theorem \ref{mainthm} follows from the lemma by the same exact reasoning used to prove Proposition \ref{hodgegen}; we will not repeat the argument here. An immediate consequence of Theorem \ref{mainthm} and Proposition \ref{hodgegen} is a formula for $\lambda$ in terms of the boundary classes and $\zeta$ or $\zeta'$.

\begin{cor}\label{hodgezeta}
In $\Pic(\hypbar_g)$,
$$
\lambda=\frac{g}{2}\zeta-\sum_{i=1}^{\lfloor\frac{g-1}{2}\rfloor}
\frac{i(i+1)}{2}\xi_i-\sum_{j=1}^{\lfloor\frac{g}{2}\rfloor}
j^2\delta_j
$$
when $g$ is even, and
$$
\lambda=g\zeta'-\sum_{i=1}^{\lfloor\frac{g-1}{2}\rfloor}
\frac{i(i+1)}{2}\xi_i-\sum_{j=1}^{\lfloor\frac{g}{2}\rfloor}
j^2\delta_j
$$
when $g$ is odd.
\end{cor}
To prove the formulas, subtract $g/2$ times (\ref{zetaeven}) from (\ref{halfformula}) for even $g$, and $g$ times (\ref{zetaodd}) from (\ref{oldformula}) for odd $g$ to obtain, respectively, $4g+2$ times the first identity, or $8g+4$ times the second one.
Since $\Pic(\hypbar_g)$ has no torsion, the result follows. Another immediate consequence of Theorem \ref{mainthm} and Corollary \ref{hodgezeta} is the following result of Gorchinskiy and Viviani; to make contact with their notation just observe that they write $\mathcal G$ to indicate our line bundle $\calz$ when $g$ is even, and our $\calz'$ when $g$ is odd.

\begin{cor}[\cite{gorvi}]\label{hodgezeta2}
When $g$ is odd, $\lambda=g\zeta'$ in $\Pic(\hyp_g)$, and $\Pic(\hyp_g)$ is generated by $\zeta'$.
When $g$ is even, $\lambda=\frac{g}{2}\zeta$ in $\Pic(\hyp_g)$, and $\Pic(\hyp_g)$ is generated by $\zeta$. If $g$ is not divisible by $4$, $\Pic(\hyp_g)$ is generated by $\lambda$ while, if $g$ is divisible by $4$, $\lambda$ generates an index two subgroup of $\Pic(\hyp_g)$.
\end{cor}
Finally, when $g$ is divisible by $4$, the greatest common divisor of the coefficients of (\ref{oldformula}) is exactly $4$. Since $\Pic(\hypbar_g)$ has no torsion, we obtain a valid identity if we divide all coefficients of (\ref{oldformula}) by $4$. This settles the question of the structure of the subgroup of $\Pic(\hypbar_g)$ generated by $\lambda$ and by the boundary classes.

\begin{prop}\label{hodgegen2}
Let $g\ge 2$ be an integer which is divisible by $4$. Then $\lambda$ and the boundary classes generate an index $2$ subgroup of $\Pic(\hypbar_g)$, and the relations between them are generated by
$$
(2g+1)\lambda=\frac{g}{4}\xi_{irr}+\sum_{i=1}^{\lfloor\frac{g-1}{2}\rfloor}
\frac{(i+1)(g-i)}{2}\xi_i+\sum_{j=1}^{\lfloor\frac{g}{2}\rfloor}
j(g-j)\delta_j
$$
\end{prop}
At this point, all we have to do to finish up is prove Lemma \ref{mainlemma}.

\medskip
\noindent \textit{Proof of Lemma \ref{mainlemma}}. The argument is essentially the one originally used to prove (\ref{oldformula}), and we shall freely employ results from \cite{corhar}. If we pick $g$ families of stable hyperelliptic curves, all with smooth and complete one-dimensional base, we can construct the $g\times g$ matrix whose entries are the degrees of the various boundary classes on the given families. It was shown in \cite{corhar} that the families can be chosen in such a way that this matrix is non-singular. This proves that the boundary classes are independent. It also reduces the task of proving (\ref{zetaeven}) modulo torsion to the one of evaluating, for each one of the $g$ families, the degrees of $\zeta$ and of the boundary classes, and showing that one gets an identity if, in (\ref{zetaeven}), one replaces each class with its degree. In other words, we must show that
\begin{equation}\label{tobeproved}
(4g+2)\deg\zeta=\deg\xi_{irr}+2\sum_{i=1}^{\lfloor\frac{g-1}{2}\rfloor}
(i+1)(2i+1)\deg\xi_i
+4\sum_{j=1}^{\lfloor\frac{g}{2}\rfloor}j(2j+1)\deg\delta_j
\end{equation}
for each one of the families.

The first family that one considers is constructed as follows. Let $D$ be a general divisor of type $(2g+2,2)$ in $Y=\PP^1\times \PP^1$, and let $f:Y\to \PP^1$ be the projection to the second factor. As $D$ is general, it is smooth, and $f_{|D}:D\to\PP^1$ is a simple covering, in the sense that above each point of $\PP^1$ there is at most one ramification point, and the ramification index at this point is $2$. Since $\calo(D)$ is a square, there is a double covering $\eta:X\to Y$ ramified at $D$. The surface $X$ is smooth and, writing $\pi$ for the composition of $\eta$ and $f$, $\pi:X\to \PP^1$ is a family of genus $g$ stable hyperelliptic curves. All singular fibers of $\pi$ are of type $\Xi_{irr}$, and the degree of $\xi_{irr}$ is equal to their number, that is, to the number of ramification points of $f_{|D}$. This can be easily calculated using the Riemann-Hurwitz formula for $f_{|D}$ and the genus formula for $D\subset Y$, and turns out to be $8g+4$. To prove (\ref{tobeproved}) in our case we must therefore show that $\deg\zeta=2$. This is also easy to do. Write $W$ for the ramification divisor of $\eta$. Since $\omega_{\pi}=\eta^*\omega_f(W)$ and $\eta^*(D)=2W$, we have that $\omega_{\pi}^{g+1}(-(g-1)W)=\eta^*\omega_f^{g+1}(2W)=\eta^*(\omega_f^{g+1}(D)) =\eta^*(\calo_{\PP^1\times \PP^1}(0,2))=\pi^*(\calo_{\PP^1}(2))$. Hence $\pi_*(\omega_{\pi}^{g+1}(-(g-1)W))=\calo_{\PP^1}(2)$, and $\deg\zeta=2$, as desired.

The remaining $g-1$ families are all obtained by the same general procedure. We start with a family of stable $(2g+2)$-pointed curves of genus zero, consisting of a family of curves $f:Y\to S$ plus sections $D_1,\dots,D_{2g+2}$, and we set $D=\sum D_i$. We assume that $S$ is a smooth complete curve, that the general fiber of $f$ is smooth, and that $\calo(D)$ is a square. We let $\eta:X\to Y$ be the double covering branched along $D$, write $R$ for the ramification divisor of $\eta$, and set $\pi=f\circ\eta$. Then $\pi:X\to S$ is a family of semistable hyperelliptic curves of genus $g$, and to obtain a family of stable hyperelliptic curves we pass to its stable model $\pi':X'\to S$. In practice, as all the degrees we need to consider are readily computed on $\pi:X\to S$, we will work mostly with this family, rather than with $\pi':X'\to S$. The degrees of the boundary classes have been computed in \cite{corhar}; we do not need to know them individually, but just that, as a consequence of Lemma (4.8) and formula (4.10) in \cite{corhar}, they are tied to the self-intersection of $D$ by the relation
\begin{align}\label{basic}
(2g+1)(D\cdot D)&=-2g\deg\xi_{irr}-\sum_{i>0}(2i+2)(2g-2i)\deg\xi_i\\
&\phantom{=}-2\sum_{i>0}(2i+1)(2g-2i+1)\deg\delta_i\notag
\end{align}
Next, we compute the degree of $\zeta$ on the family $\pi':X'\to S$, that is, $\deg(\calz_{\pi'})$. The line bundle $\calz_{\pi'}$ is the pushforward of the line bundle on $X$
$$
L=\omega_\pi^{g+1}\left(-(g-1)R-\sum(2g-4i)G_i-\sum(g-2i-1)E_i\right)
$$
where $E_i$, $0\le i\le \lfloor(g-1)/2\rfloor$, is the sum of all tails of type $\Xi_i$ in the fibers of $\pi$, each counted with the appropriate multiplicity, and $G_i$ is the sum of all tails of type $\Delta_i$, also counted with the appropriate multiplicity. By this we mean the following. Let $p$ be a node in a singular fiber of $\pi$; complex analytically, $X$ is of the form $xy=t^m$ near $p$, where $t$ is a local parameter on $S$. Either $p$ is a node of type $\Delta_i$, or it belongs to a pair of type $\Xi_i$; call $T$ the corresponding tail. Then $T$ appears with multiplicity $m$ in $G_i$ or in $E_i$. Notice that the contribution of node $p$ to the degree of $\delta_i$, or the contribution to the degree of $\xi_i$, $i\ge 1$, of the pair of nodes to which $p$ belongs, is $m$; instead, when $i=0$, we get a contribution of $2m$ to the degree of $\xi_{irr}$. Notice also that the intersection number of $R$ with $T$ is $2i+1$ if $p$ is of type $\Delta_i$, and $2i+2$ if $p$ belongs to a pair of type $\Xi_i$. The appearance of tails of type $\Xi_0$ in the definition of $L$ is due to the fact that the pullback to $X$ of the Weierstrass divisor $W_{\pi'}$ in $X'$ is $R+E_0$. We also observe that $(D\cdot D)=-(\omega_f\cdot D)$, as $D$ is a disjoint union of sections. In view of these remarks, and since the degree of $\calz_{\pi'}=\pi_*(L)$ is equal to the intersection number of $L$ with a section of $\pi$, we have that
\begin{align}
(2g+2)\deg\zeta &=(L\cdot R)=(\omega_f^{g+1}\cdot D)+2(R\cdot R)-(g-1)(E_0\cdot R)\notag\\
&\phantom{=}-\sum(g-2i-1)(E_i\cdot R)-\sum(2g-4i)(G_i\cdot R)\notag\\
&=-g(D\cdot D)-(g-1)\deg\xi_{irr}-\sum(g-2i-1)(2i+2)\deg\xi_i\notag\\
&\phantom{=}-\sum(2g-4i)(2i+1)\deg\delta_i\notag
\end{align}
Identity (\ref{tobeproved}) follows by substituting for $(D\cdot D)$ the value given by (\ref{basic}).

\hfill q.e.d.

\end{document}